\documentclass{article}
\usepackage[utf8]{inputenc}
\usepackage[T2A]{fontenc}
\usepackage{amsmath}
\title{AVOIDABLE WORDS}
\date{}
\author{IRINA MELNICHUK}
\begin{document}

\maketitle
\setlength{\parindent}{30ex}                       ABSTRACT \\  
\setlength{\parindent}{20ex}The set of all avoidable patterns in n \setlength{\parindent}{20ex}or   fewer letters can be avoided on an \setlength{\parindent}{20ex}alphabet with 2(n+2) letters. \\ \\ 

\setlength{\parindent}{5ex} Слово \textit{а} называется блокирующим [1], если для любого конечного алфавита множество слов в этом алфавите не содержащих значений слова \textit{а}, конечно. Напомним, что значением слова \textit{а} называется результат подстановки некоторых слов вместо букв слова \textit{а}, при этом одинаковые буквы заменяются одинаковыми словами. Неблокирующее слово  \textit{а} называется n-избегаемым, если в некотором n-буквенном алфавите (а, значит, и в каждом n-буквенном ) существует бесконечное множество слов, не содержащих значений слова \textit{а}. В проблеме 3 из [2] ставится вопрос об определении для произвольного неблокирующего слова наименьшего алфавита, в котором это слово избегаемо. Известны некоторые частные результаты, относящиеся к этому вопросу.
 
   В [3] показано, что слово \(x^2\) избегаемо в 3-буквенном алфавите, а \(x^3\) избегаемо в 2-буквенном алфавите. В [4] установлено, что множество всех двойных слов над n-буквеннм алфвитом  избегаемо в алфавите, содержащем 3[n/2]+3 букв,где [n/2]-целая часть n/2.  В [5] построена последовательность слов в алфавите из четырех букв, которая избегает любое полное слово , и показано, что в алфавите из из трех букв такой последовательности не существует. В [6] в теоремах 1.2 и 1.3 дана верхняя  граница мощности алфавита, в котором избегаемо произвольное неблокирующее слово.
   
     Пусть \textit{а} - неблокирующее слово, \(\alpha=\alpha(a)\) - число различных букв, входящих в \textit{а}. Обозначим через \textit{m}=\textit{m}(a) минимальное число \textit{n} такое, что \textit{а} является n-избегаемым. По теореме 1.2 из [6] для любого неблокирующего слова  \textit{а}, 
     
          \(m<4*(\alpha+2)*[ln(\alpha+2)]\). \\ 
     По теореме 1.3 из[6],
      
          \( m<9*(\alpha+20)\) для любого неблокирующего слова \textit{а}. \\   
     Цель статьи - следующее улучшение  оценки числа m.
      
      ТЕОРЕМА.  Для любого неблокирующего слова \textit{u}  \\  
          \( m(u)<=2*\alpha(u)+4 \)  \\ 
      
      Д О К А З А Т Е Л Ь С Т В О.  Обозначим через \textit{k} то из чисел  \(2*\alpha(u)+2\),  \(2*\alpha(u)+4\), которое делится на 4. докажем, что \textit{u} является k-избегаемым. Будем предполагать, что \(\alpha(u)>1\), так как для \(\alpha(u)=1\) утверждение теоремы вытекает из [3].
      
       Нам понадобится описание блокирующих слов, полученное в [1]. Введем в рассмотрение слова над алфавитом \(\Xi=\big\{\xi_1,\xi_2,...\big\}\). Полагаем \(Z_1=\xi_1\) и далее по индукции \(Z_{n+1}=Z_n\xi_{n+1}Z_n\), где n=1,2,.... Обозначим через \(F_1(A)\) - множество всех непустых слов над алфавитом A, через  \(\alpha(a)\) - множество всех букв слова \textit{а}. 
       
       Для слова \textit{а} отображение  \(f:\alpha(\textit{а})\rightarrow F_1(\Xi)\) такое, что f(\textit{а}) является подсловом \(Z_n\) для некоторого n, назовем В-отображение слова \textit{а}. А.И. Зимин в [1] доказал, что слово \textit{а} является блокирующим тогда и только тогда, когда существует В-отображение этого слова.
       
       Рассмотрим в симметрической группе \(S_k\)  перестановки \\ 
\begin{align*}
f_1 = (1,3,5,...,k-1),\quad g_1 = (2,6,...,k-2)   \\ 
\setlength{\parindent}{9ex}    f_2 = \setlength{\parindent}{9ex}f_1^2,\qquad\qquad\qquad\qquad\quad g_2 \setlength{\parindent}{9ex}= g_1 f_1  \\
\setlength{\parindent}{9ex}    ...\qquad\qquad\qquad\qquad\qquad\qquad\quad... \\ 
     f_\frac{k}{2} = f_1^\frac{k}{2},\qquad\qquad\qquad\qquad\quad g_\frac{k}{2} = g_1 f_1^{\frac{k}{2}-1} \\ 
    \end{align*}
    
    Получим k различных перестановок, обозначим их через 
 \(v_1,...,v_k\). Для выбранных перестановок зафиксируем  k слов \(a_1,...,a_k\) над алфавитом \(X = \big\{x_1,...,x_k\big\}\)  полагая \(a_i = x_{v_i(1)} ... x_{v_i(k)}\). \\ 
 О П Р Е Д Е Л Е Н И Е 1. Зададим отображение  \(\varphi_0:X\rightarrow F_1(X)\) по правилу   \(\varphi_0(x_i)=a_i\). Обозначим  \(J_1 = \varphi_0(x_1)\). Пусть  \(J_{m-1}\) уже определено, тогда полагаем \(J_m = \varphi_0(J_{m-1})\). \\  
 Для доказательства теоремы достаточно показать,что для любого  \(m\geq1\) слово \(J_m\)не содержит значений неблокирующего слова  \(u\). Доказательство этого факта вытекает из следующего ниже предложения. \\ 
П Р Е Д Л О Ж Е Н И Е. Пусть при отображении  \(\varphi_0:\alpha(u)\rightarrow F_1(X)\) для слова  \(u\)  содержащего не более  \(\frac{k}{2}-1\) букв, слово \(\varphi(u)\)  является подсловом \(J_n\) для некоторого \(n\geq1\). Тогда существует  B-отображение \(u\), т.е. \(u\)  слово блокирующее \\
 Д О К А З А Т Е Л Ь С Т В О. \\
 Будем рассматривать множество \(A = \big\{a_1,...,a_k\big\}\) как некоторый алфавит. Тогда  для  \(n\geq1\) слово \(J_n\) можно рассматривать как слово над алфавитом  A, назовем его в этом случае цепью. Вхождение  \(b\) из \(F_1(A)\) в цепь  \(J_n\) подцепью  \(b\), в частности, вхождение \(a_i\) в цепь назовем звеном \(a_i\). \\
Таким образом, подцепь  \(a_i\) это вхождение слова из \(F_1(A)\)  длины 1 в \(J_n\), слово \(a_i\) слово из  \(F_1(X)\) длины  k. Будем говорить, что  слово  \(a\) из  \(F_1(X)\) является пословом подцепи \(b\),если соотвествующее подцепи  \(b\) слово в алфавите X содержит подслово \(a\). Будем говорить. что подцепи  \(C\) и \(D\) графически равны (\(C\equiv D\))  если они являются вхождениями графически равных слов. \\  
  Пусть дано некоторое вхождение слова \textit{а} в \(J_n\). Замыканием вхождения  \(a\) назовем наименьшую подцепь цепи   \(J_n\), содержащую данное вхождение  \(a\).  Обозначим  замыкание  \(a\) через  \( \big[a\big] \). \\ 
     Доказательство предложения проведем индукцией по длине  
     \(\big[\varphi(u)\big]\).\\ 
   Ясно, что если  \(\big[\varphi(u)\big]\) совпадает с одним из звеньев \(a_i, \: i\in\big\{1,...,k\}\), то \(u\) - блокирующее слово. \\
     Индукционный шаг. Предположим, что \(\varphi(u)\)  является подсловом \(J_n\) для некоторого \(n\geq1\)  и подцепь  \(\big[\varphi(u)\big]\)содержит  более одного звена. Докажем, используя индукционное предположение, что существует B-отображение слова \(u\).  \\ 
     О П Р Е Д Е Л Е Н И Е  2. Слово \(b\)  из  \(F_1(X)\) назовем базисным, если оно является подсловом  \(a_i\)  для некоторого \(i\in\big\{1,...,k\big\}\)  и содержит  не менее двух букв  x  с четными индексами \\ 
  О П Р Е Д Е Л Е Н И Е  3  . Если \(a_i = b_1 b_2\),  \(a_j = c_1 c_2\), где  \(b_2, c_1 \in F_1(X)\), то слово \(b_2 c_1\) назовем смежным. \\ 
   О П Р Е Д Е Л Е Н И Е  4. Подслово b = cd, \(c \neq 	\wedge , d \neq 	\wedge \) слова  \(a_i, i\in\big\{1,...,k\big\}\) назовем разбитым, если для некоторого  вхождения  буквы \(y\) в \(u\)  соответствующее  замыкание \(\varphi(y)\) содержит последним (первым) звеном \(a_i\), причем   c (d) является правым (левым) подсловом \(\varphi(y)\). \\ 
 О П Р Е Д Е Л Е Н И Е    5  . Для каждого  r-го вхождения буквы \(y\) в слово  \(u\) зафиксируем подцепь  \(\big[\varphi(y)\big]_r\),
 являющуюся замыканием r-го  вхождения   \(\varphi(y)\) в \(J_n\), соответствущего  r-му вхождению \(y\) в \(u\). Предположим, что  \(\big[\varphi(y)\big]_r = c\varphi(y) d\). Если \(c \neq 	\wedge \), то первое звено  \(\big[\varphi(y)\big]_r\) назовем разбитым. Если \(d \neq 	\wedge \), то последнее звено назовем разбитым, другие  звенья \(\big[\varphi(y)\big]_r\)  назовем неразбитыми.
 Для краткости \(\big[\varphi(y)\big]_1\) обозначим \(\big[\varphi(y)\big]\). \\ 
О П Р Е Д Е Л Е Н И Е   6. Букву \(x\) с четным индексом слова \(a_i\) назовем базисной, если не существует такого вхождения  буквы \(y\) слова  \(u\), что звено \(a_i\) является первым для соответствующего замыкания \(\varphi(y)\) и слово  \(\varphi(y)\)  содержит \(x\) слева первой или второй буквой. \\ 
  Из определения вытекает, что если  \(x\) - базисная буква слова  \(a_i\), и звено \(a_i\) является первым в цепи \(\big[\varphi(y)\big]_r\), то либо вхождение базисного слова вида  \(x_i x_j x\) в  \(a_i\)  является вхождением в  \(\varphi(y)\), либо вхождение  \(x\) в  \(a_i\)  не является вхождением в \(\varphi(y)\). \\ 
Л Е ЕМ М А  1. Слова  \(a_1, ..., a_k\) обладают следующими свойствами: \\ 
a) каждое базисное слово является подсловом \(a_i\) для  единственного  \(i\) из множества  \(\big\{1, ..., k\big\}\); \\ 
б) если \(b\) - левое подслово \(a_i\) , \(c\) левое подслово  \(a_j\) и \(b = c, \; |b| \geq 2\), то \(a_i = a_j\);\\
в) не существует смежного слова, которое является подсловом \(a_i\) для некоторого  \(i\in\big\{1,...,k\big\}\); \\ 
г) в слове  \(J_n\)  нет одинаковых рядом стоящих подслов. \\ 
д) у каждого слова  \(a_i\) есть базисная буква. \\ 
  Доказательство пунктов а)-г) легко вытекает из определения слов 
  \(a_i\). Для доказательства утверждения д) заметим, что число четных букв в слове  \(a_i\) не  меньше  \(\alpha(u)+1\), а число букв слова \(u\) равно \(\alpha(u)\). \\ 
О П Р Е Д Е Л Е Н И Е  7. Для каждого слова \(a_i,\; i\in\big\{1,...,k\big\}\) зафиксируем базисную букву  \(a_i\) и назовем ее основной базисной буквой   \(a_i\). \\ 
О П Р Е Д Е Л Е Н И Е   8. Зафиксируем отображение  \(C: F_1(X)\rightarrow F_1(\Xi)\) следующим образом. Представим  натуральное  число \(p\) в виде  \(p=2^i+r2^{i+1}\)  для некоторых целых \(r\geq 0, i\geq 0\). Такое представление однозначно. Полагаем  \(C(x_p)=\xi_{i+1}\) для любого  \(x_p\), \\ где \(1\leq p\leq k-1, \; p=2^i+r2^{i+1}\). Тогда слово \(C(x_1 x_2...x_{k-1})\)совпадает для некоторого числа \(t\)  с левым подсловом  \(Z_t\) . Выберем наименьшее такое  \(t\). Пусть  \(Z_t=C(x_1 x_2...x_{k-1})Z'\), тогда полагаем  \(C(x_k)=Z'\xi_{t+1}\). \\  
З А М Е Ч А Н И Е. По определению,  \\  
\(C(x_1)=C(x_3)=...=C(x_{k-1})=\xi_1\),   
\(C(x_2)=C(x_6)=...=\xi_2\), \\ 
\(C(x_1 x_2...x_k)=Z_t\xi_{t+1}\). \\ 
Так как слово \(Z_t\) от перестановки между собой букв \(\xi_1\) или   \(\xi_2\)  не меняется, то для любого \(a_i, \: i\in\big\{1,...,k\big\}\),  \(C(a_i)=Z_t\xi_{t+1}\). \\  
Л Е М М А  2. Пусть для каждого  \(y\) из  \(\alpha(u)\) слово  \(\varphi(y)\) удовлетворяет условию : если \(a_i\)  - подслово   \(\big[\varphi(y)\big]\), то вхождение основной базисной буквы   \(a_i\) не является вхождением в  \(\varphi(y)\).  Тогда  \(u\) - блокирующее. \\  
Д О К А З А Т Е Л Ь С Т В О. Пусть отображение \(\varphi\) удовлетворяет условию леммы. В этом случае  каждое из слов  \(a_i, \: i\in\big\{1,...,k\big\}\) не является подсловом  \(\varphi(y)\), так как по лемме 1 д), у каждого \(a_i\)  есть основная базисная буква. Значит, \(\varphi(y)\) является для некоторых i, j подсловом  \(a_ia_j\), причем по крайней мере последняя буква   \(a_j\)-  буква  \(x_k\) не входит в  \(\varphi(y)\). Поэтому   \(C(\varphi(y))\)  является подсловом  \(Z_t\xi_{t+1}Z_t\), то есть является блокирующим. Лемма 2 доказана. \\  
   Будем предполагать в дальнейшем, что условие леммы  2 не выполняется, то есть для некоторого   \(y\) из  \(\alpha(u)\)  слово \(\varphi(y)\)  содержит вхождение  основной базисной буквы некоторого  \(a_i\) из \(\big[\varphi(u)\big]\).  Из определения  6 основной базисной буквы слова следует, что если для некоторого p-го вхождения  \(y\) в  \(u\)  цепь  \(\big[\varphi(y)\big]_p\) начинается (заканчивается) звеном  \(a_i\) , причем вхождение основной базисной буквы слова  \(a_i\) является вхождением  \(\varphi(y)\),  то для любого r-го вхождения  \(y\) в \(u\) цепь \(\big[\varphi(y)\big]_r\)  начинается (заканчивается) подсловом   \(a_i\). \\ 
 О П Р Е Д Е Л Е Н И Е   9. Слову \(u\)  поставим в соответствие  слово \(u_1\) следующим образом: вычеркнем из   \(u\) вхождения  всех таких букв \(y\)  для которых \(\big[\varphi(u)\big]\) состоит из разбитых звеньев (и, значит , содержит одно или два звена), причем вхождение основной базисной буквы каждого из этих звеньев не является вхождением в  \(\varphi(y)\). \\ 
   На множестве  P  всех звеньев цепи \(\big[\varphi(u)\big]\), построим функцию  \\ 
\(f:P\rightarrow \big\{R,L,0\big\}\). \\ 
 Предположим, что для некоторого  r-го вхождения \(x\) в \(u\) подцепь  \(\big[\varphi(x)\big]_r\) содержит звено  \(a_i\).  \\ 
  а)Если \(a_i\) неразбитое (определение 5), или  \(\big[\varphi(x)\big]_r\)=\(a_i\)  для  \(x\) из \(\alpha(u_1)\), то полагаем  \(f(a_i)=0\).  \\  
 б)  разбитое звено  \(a_i\) - первое из подцепи \(\big[\varphi(x)\big]_r\), \(\varphi(x)=sb\)  где s - правое подслово \(a_i\). Если основная базисная буква  \(a_i\) входит в s, то полагаем \(f(a_i)=R\); \\ 
 в) разбитое  звено \(a_i\) - последнее в подцепи  \(\big[\varphi(x)\big]_r\), \(\varphi(x)=bs\) где s - левое подслово \(a_i\). Если основная базисная буква \(a_i\) входит в s, то полагаем \(f(a_i)=L\); \\ 
 г) пусть для разбитого звена \(a_i\)  - не выпоняются условия пунктов а)-в), то есть не определена функция f, тогда  полагаем   \(f(a_i)=0\). Очевидно, что этот случай возможен, когда  \(a_i\) - крайнее звено в цепи  \(\big[\varphi(u)\big]\) и основная базисная буква \(a_i\) не входит в  \(\varphi(u)\).  \\ 
    Для каждого разбитого звена существует единственная основная базисная буква, поэтому определение функции корректно. \\ 
    З А М Е Ч А Н И Е. Для любого r-го вхождения  \(y\) в  \(u_1\) значение f для последнего (первого) звена \(\big[\varphi(y)\big]_r\) равно значению f для последнего (первого) звена\(\big[\varphi(y)\big]\). \\ 
    Действительно, предположим,   \(y\in\alpha(u_1)\), 
    \(\big[\varphi(y)\big]_r\)=b\(a_j\),
 \(\big[\varphi(y)\big]\)=d\(a_i\). \\ 
  Пусть  \(f(a_i)=L\), \(a_i=a_{i_1} a_{i_2}\),  где  \(a_{i_1}\) правое подслово  \(\varphi(y)\).  Тогда, по определению f, основная базисная буква \(a_i\)  входит в  \(a_{i_1}\). Так как базисная буква имеет четный индекс, то длина \(a_{i_1}\) больше или равна двум и по лемме 1 б) \(a_i=a_j\)  и  \(f(a_i)=L=f(a_j)=L\) . \\  
  Пусть \(f(a_i)=R\),  \(\big[\varphi(y)\big]_r=a_j\)b,  
 \(\big[\varphi(y)\big]=a_i\)d,  \(a_i=a_{i_1} a_{i_2}\),  где  \(a_{i_2}\) есть левое подслово  \(\varphi(y)\). Тогда основная базисная буква \(a_i\)  входит в \(a_{i_2}\), и  значит, \(a_{i_2}\) содержит базисное подслово, однозначно определяющее  слово   \(a_i\) и \(a_i=a_j\)  . \\  
 Если   \(f(a_i)=0\) и  \(a_i\) неразбитое (определение 5), или  \(\big[\varphi(y)\big]_r\)=\(a_i\)  для  \(y\) из \(\alpha(u_1)\), то \(\big[\varphi(y)\big]\)  содержит базисное подслово \(a_i\), и значит, для любого r-го вхождения y в \(u_1\), \(\big[\varphi(y)\big]_r\)=\(a_i\). 

  Построим отображение  \(\Psi:\alpha(u_1)\rightarrow F_1(A)\). Пусть  \(y\in\alpha(u_1)\).  \\  
  a) Предположим, \(\big[\varphi(y)\big]=a_i\). Полагаем   \(\Psi(y)=a_i\). \\ 
  b) Предположим,  подцепь \(\big[\varphi(y)\big]\) \ содержит более одного звена, первое звено этой подцепи  \(a_i\), последнее  \(a_j\), то есть , \(\big[\varphi(y)\big]=a_iBa_j\). Полагаем   \\   \\

\(\Psi(y) =
  \begin{cases}
  a_i B a_j, \;if\; f(a_i)=R, f(a_j)=L, \\
  a_i B, \;if\; f(a_i)=R, f(a_j)\neq L,\\
  B a_j, \;if\; f(a_i)\neq R, f(a_j)=L,\\
  B, \;if\; f(a_i)\neq R, f(a_j)\neq L.
  \end{cases}
\) \\ 
   Так как каждое разбитое или неразбитое звено из подцепи \(\big[\varphi(y)\big]\),кроме. может быть, первого и последнего, входит в  \(\big[\Psi(y)\big]_r\)  для единственных  \(r \geq 1\),  \(y\)  из  \(\alpha(u_1)\) ,то  \(\Psi(u_1)\) - подслово  \(J_n\) и по определению отображения  \(\Psi\), для каждого   \(y\) из \(\alpha(u_1)\), \(\Psi(y)= \big[\Psi(y)\big]\). \\ 
    По отбражению   \(\Psi\) , построим отображение   \(\Psi_1:\alpha(u_1)\rightarrow F_1(A)\) , полагая, если  \(\Psi(y)=a_{i_1}...a_{i_m}\) , то  \(\Psi_1(y)=x_{i_1}...x_{i_m}\)  \\    Тогда \(\Psi_1(y)\) - подслово  \(J_{n-1}\) причем  длина  \(\big[\Psi_1(u_1)\big]\) меньше длины \(\big[\Psi(u)\big]\),  поэтому, по индукционному предположению, \(u_1\)  - слово блокирующее и существует  B-отображение слова  \(u_1\). \\ 
    Пусть  \(G(u_1)\)  - подслово  \(Z_r\). Построим отображение \(G_1:\alpha(u)\rightarrow F_1(\Xi)\)  и докажем, что это B-отображение слова \(u_1\). \\ 
    Введем отображение  \(H:\Xi\rightarrow F_1(\Xi)\)  , полагая \(H(\xi_i)=Z_t\xi_{i+t+1}Z_t\xi_{t+1}\) (t взято из определения   8). Легко видеть, что  H(\(Z_r\)) подслово  \(Z_{r+t+1}\) и значит, H(G(\(u_1)\)) - подслово  \(Z_{r+t+1}\). \\ 
     О П Р Е Д Е Л Е Н И Е    10. Для каждого  \(y\notin \alpha(u_1)\)  полагаем  \(G_1(y)=C(\varphi(y))\) (отображение C дано в определении  8). \\ 
     Пусть  \(y\in \alpha(u_1)\).  \\ 
 a) \(\varphi(y)\) - подслово \(a_i\)  для некоторого   \(i\in\big\{1,...,k\big\}\).  Предположим, \(a_i = m\varphi(y)s\),   H(G(\(y)\))=\(Z_tbZ_t\xi_{t+1}\).   Тогда полагаем,  \(G_1(y)=pbC(m\varphi(y))\), где через p обозначено слово, полученное из слова C(\(\varphi(y)s)\) отбрасыванием последней буквы  \(\xi_{t+1}\). \\ 
 Заметим, что H(G(y)) - подслово  \(Z_{r+t+1}\) и \(G_1(y)\)  - подслово  H(G(y)). \\ 
 b) Пусть цепь \(\big[\varphi(y)\big]\)  содержит более одного звена, первое звено этой подцепи    \(a_i\), последнее  \(a_j\), то есть  \(\big[\varphi(y)\big]\)=\(a_iBa_j\), \(\varphi(y)\)=mBs.  Пусть  H(G(\(y)\))=\(Z_tbZ_t\xi_{t+1}\). Обозначим через  p слово, полученное из C(m) отбрасыванием последней буквы \(\xi_{t+1}\). Полагаем :  \\   \\

\(G_1(x) =
  \begin{cases}
pbC(s),                        \;if\; f(a_i)=R, f(a_j)=L, \\
pbZ_t\xi_{t+1}C(s),         \;if\; f(a_i)=R, f(a_j)\neq L,\\
C(m)Z_tbC(s),              \;if\; f(a_i)\neq R, f(a_j)=L,\\
C(m)Z_tbZ_t\xi_{t+1}C(s),  \;if\; f(a_i)\neq R, f(a_j)\neq L.
  \end{cases}
\) \\ 
Л Е М М А  3. Для построенного отбражения  \(G_1:\alpha(u)\rightarrow F_1(\Xi)\), слово  \(G_1(u)\) является  подсловом   \(Z_t\xi_{t+1}H(G(u_1))Z_t\). \\  
Д О К А З А Т Е Л Ь С Т В О леммы  проведем  индукцией  по длине слова \(u_1\). По лемме 2 слово  \(u_1\) не пустое. Пусть \(u_1\)- однобуквенное слово  x, H(G(x))=\(Z_tb_tZ_t\xi_{t+1}\), и \(u_1\) получено из слова  u=dxh  вычеркиванием подслов  d, h. Каждая из цепей  \(\big[\varphi(d)\big]\),  \(\big[\varphi(h)\big]\)  содержит не более двух звеньев. Для определенности, пусть \(\big[\varphi(d)\big]\)=\(a_i\)  и  \(\big[\varphi(h)\big]\)=\(a_ja_m\) (остальные случаи рассматриваются аналогично) . Запишим \(a_i\), \(a_ja_m\) в виде 
\(a_i\)=\(k_1\varphi(d)k_2\), \(a_ja_m\)=\(r_1\varphi(h)r_2\).   
  Так как буквы подслова h были вычеркнуты из слова u, то \(\varphi(h)\) не содержит вхождения основной базисной буквы   \(a_j\) и не содержит вхождения основной базисной буквы \(a_m\). Предположим,  \(f(a_i)\neq R\)(аналогично рассматривается случай, когда  \(f(a_i)=R\) ).  
  Так как звено  \(a_j\)  не последнее в цепи   \(\big[\varphi(u)\big]\)  и буквы слова h  вычеркиваются из u, то \(f(a_j)=L\) и \(f(a_m)=0\). \\  
  По определению  10,   \(G_1(x)\)=C(\(k_2)Z_tbC(r_1)\), \(G_1(d)=C(\varphi(d))\), \(G_1(h)=C(\varphi(h))\) \\  
   Поэтому \(G_1(dxh)=C(\varphi(d))C(k_2)Z_tbC(r_1)C(\varphi(h))\),  где \(C(\varphi(d))C(k_2)\) - правое подслово \(C(a_i\))=\(Z_t\xi_{t+1}\), а  \(C(r_1)C(\varphi(h))\) - левое подслово   \(Z_t\xi_{t+1}Z_t\),  содержащее  \(Z_t\xi_{t+1}\). \\  
   Базис индукции доказан. \\  
   Предположим, слово \(u_1\) имеет вид  \(u_1=s_1xyw_1\),где    \(x, y\in\alpha(u_1)\) и  \(u_1\) получено из слова u=sxvyw вычеркиванием послова  v и некоторых других букв из слов   s, w. Пусть \(G(u_1)=b_1b_2\) - подслово  \(Z_r\), где  \(G(s_1x)=b_1\), \(G(yw_1)=b_2\). По индукционному предположению, \(G_1(sx)\) - подслово  \(Z_t\xi_{t+1}H(G(s_1x))Z_t\) и  \(G_1(yw)\)  - подслово   \(Z_t\xi_{t+1}H(G(yw_1))Z_t\). Буквы  слова v вычеркиваются из u при получении  \(u_1\), поэтому   \(\big[\varphi(v)\big]\) содержит не более двух звеньев.  Предположим,    \(\big[\varphi(v)\big]\)=\(a_i\) (случай, когда   \(\big[\varphi(v)\big]\) пустое или содержит два звена, рассматривается аналогично).  Запишем слово   \(a_i\) в виде: \(a_i=k_1\big[\varphi(v)\big]k_2\). Пусть \(H(G(s_1x))=Z_tb_1Z_t\xi_{t+1}\), \(H(G(yw_1))=Z_tb_2Z_t\xi_{t+1}\). \\ 
a) Предположим, \(f(a_i)=R\). Тогда по определению  10,  \(G_1(sx)=d_1b_1Z_t\xi_{t+1}C(k_1)\), \(G_1(yw)= pb_2d_2\) где  p - слово полученное из  \(C(k_2)\) (определение  7) отбрасыванием последней буквы \(\xi_{t+1}\) и \(d_1\),\(d_2\) - подслова  \(Z_{t+1}\).  Но  \(C(k_1)C(\varphi(v))p=Z_t\),   и значит,   \(G_1(sxvyw)=d_1b_1Z_{t+1}b_2d_2\). \\ 
б) Предположим, \(f(a_i)=L\). Тогда, по определению  10, \(G_1(sx)=d_1b_1C(k_1)\), \(G_1(yw)=C(k_2)Z_tb_2d_2\). Но  \(C(k_1\big[\varphi(v)\big]k_2)=Z_t\xi_{t+1}\) поэтому    \(G_1(sxvyw)=d_1b_1Z_{t+1}b_2d_2\). \\ 
в) Предположим  \(f(a_i)=0\). Так как  \(a_i\)  разбитое звено, не являющееся крайним в цепи   \(\big[\Psi_1(u_1)\big]\) ,   то   либо  \(\big[\varphi(x)\big]=a_i\), либо  \(\big[\varphi(y)\big]=a_i\). \\  
Если   \(\big[\varphi(x)\big]=a_i\), доказательство такое же, как в случае   б).  Если   \(\big[\varphi(y)\big]=a_i\), то доказательство аналогично случаю  a). \\
 Лемма  3 доказана. \\  
Таким образом, отображение  \(G_1:\alpha(u)\rightarrow F_1(\Xi)\)  является В-отображением  слова u, и предложение доказано.  \\  
  Отсюда вытекает, что если  u - неблокирующее слово с числом различных букв  \(\alpha(u)\), то бесконечное множество слов   \(J_m\), m=1,2 ,... не содержит значений слова  u, то есть , u избегаемо в алфавите, содержащем   2\(\alpha(u)\)+4 буквы.  \\ 
 Теорема доказана. \\ \\ 
  Доказательство этой  теоремы было получено давно,в таком варианте было послано в журнал  в 1996 году, но не опубликовано. Благодарю Д.Макналти, К.Адаричеву за то, что они меня разыскали и посоветовали опубликовать этот результат.  \\ \\

              СПИСОК ЦИТИРОВАННОЙ ЛИТЕРАТУРЫ \\ \\ 
              
[1] ЗИМИН А.И.  Блокирующие множества термов.// Матем.сб., 1982. Т.119, 3, с. 363-375. \\ 

[2] BEAN D., EHRENFEUCHT A.,McNULTY G. Avoidable patterns in strings of symbols.// Pacific J. Math. 1979. v 84. N2, p.261-294. \\ 

[3] АРШОН Е.С. Доказательство существования n-значных бесконечных ассиметрических последовательностей. // Матем.сб. 1939. Т.2)44),4, С.769-779.   \\ 

[4] МЕЛЬНИЧУК И.Л. Существование бесконечных конечно порожденных свободных полугрупп в некоторых многообразиях полугрупп.// Алгебраические системы с одним действием и отношением. Ленинград: изд-во ленинградского пединститута, 1985, с. 74-83. \\ 

[5 ] ПЕТРОВ А. Н. Последовательность, которая избегает любое полное слово // Матем.заметки.1988.Т44,4. С. 517-522.  \\

[6]  BAKER K. A.,McNULTY G.F., TAYLOR W.  Growth Problems for Avoidable Words //  Teoret. Comput.Science. 1989. v.69. P. 319-345. \\  
              
\end{document}